\documentclass[12pt,oneside]{amsart}
\usepackage{amsfonts}
\usepackage{amscd}
\usepackage{amsthm}
\usepackage{amssymb}
\usepackage{amsmath}
\usepackage{epsfig}
\usepackage{graphicx}

\title{Blocking Light in Compact Riemannian Manifolds}
\author{Jean-Fran\c{c}ois Lafont}
\address{Department of Mathematics,
The Ohio State University,
Columbus, OH 43210}
\email{jlafont@math.ohio-state.edu}
\author{Benjamin Schmidt}
\address{Department of Mathematics,
University of Michigan,
2074 East Hall, 530 Church St.,
Ann Arbor, MI 48109-1043}
\email{bischmid@umich.edu}

\theoremstyle{proposition}
\newtheorem{Lem}{Lemma}[section]
\newtheorem{Prop}[Lem]{Proposition}

\newtheorem*{Def}{Definition}

\theoremstyle{plain}

\newtheorem{Thm}[Lem]{Theorem}
\newtheorem{Thm2}{Theorem}
\newtheorem{Cor}[Lem]{Corollary}
\newtheorem*{Conj}{Conjecture}

\theoremstyle{remark}

\newtheorem*{Prf}{Proof}

\DeclareMathOperator{\Diam}{Diam}
\DeclareMathOperator{\length}{length}
\DeclareMathOperator{\Cut}{Cut}
\DeclareMathOperator{\interior}{int}
\DeclareMathOperator{\vol}{vol}
\DeclareMathOperator{\inj}{inj}
\DeclareMathOperator{\Light}{Light}

\newcommand{\mA}{\mathcal A}
\newcommand{\mS}{\mathcal S}
\newcommand{\mZ}{\mathbb Z}
\newcommand{\Lm}{\Lambda}
\newcommand{\Lp}{\Lambda ^\prime}
\newcommand{\mP}{\mathcal P}
\newcommand{\mR}{\mathbb R}

\begin{document}

\begin{abstract}
We study compact Riemannian manifolds $(M,g)$ for which the light from any given point $x\in M$ can be shaded away from any other point $y\in M$ by finitely many point shades in $M$.  Compact flat Riemannian manifolds are known to have this finite blocking property.  We conjecture that amongst compact Riemannian manifolds this finite blocking property characterizes the flat metrics.  Using entropy considerations, we verify this conjecture amongst metrics with nonpositive sectional curvatures.  Using the same approach, K. Burns and E. Gutkin have independently obtained this result.  Additionally, we show that compact quotients of Euclidean buildings have the finite blocking property.

On the positive curvature side, we conjecture that compact Riemannian manifolds with the same blocking properties as compact rank one symmetric spaces are necessarily isometric to a compact rank one symmetric space.  We include some results providing evidence for this conjecture.    
\end{abstract}

\maketitle

\section{Introduction}

To what extent does the collision of light determine the global geometry of space?  In this paper we study compact Riemannian manifolds with this question in mind.  Throughout, we assume that $(M,g)$ is a smooth, connected, and compact manifold equipped with a smooth Riemannian metric $g$.  Unless stated otherwise, geodesic segments $\gamma \subset M$ will be identified with their unit speed paramaterization $\gamma : [0,L_{\gamma}] \rightarrow M,$ where $L_{\gamma}$ is the length of the segment $\gamma$.   By the \textit{interior} of a geodesic segment $\gamma$ we mean the set $\interior(\gamma):=\gamma((0,L_{\gamma}))\subset M$.

 \begin{Def}[Light]
Let $X,Y \subset (M,g)$ be two nonempty subsets, and let $G_{g}(X,Y)$ denote the set of geodesic segments  $\gamma \subset M$ with initial point $\gamma(0) \in X$ and terminal point $\gamma(L_{\gamma})\in Y$.  The \textit{light from X to Y} is the set $$L_{g}(X,Y)=\{\gamma \in G_{g}(X,Y) | \interior(\gamma) \cap (X \cup Y)= \emptyset \}.$$  
\end{Def}

\begin{Def}[Blocking Set]
Let $X,Y \subset M$ be two nonempty subsets.  A subset $B\subset M$ is a \textit{blocking set} for $L_{g}(X,Y)$ provided that for every  $\gamma \in L_{g}(X,Y),$ $$\interior(\gamma) \cap B \neq \emptyset.$$  
\end{Def}

In this paper we focus on compact Riemannian manifolds for which the light between pairs of points in $M$ is blocked by a finite set of points.  We remark that by a celebrated theorem of Serre \cite{Se}, $G_{g}(x,y)$ is infinite when $x,y \in M$ are two distinct points.  In contrast, $L_{g}(x,y) \subset G_{g}(x,y)$ may or may not be a infinite subset.  For example, in the case of a round metric on a sphere, $|L(x,y)|$ is infinite only when $x=y$ or $x$ and $y$ are an antipodal pair.  

\begin{Def}[Blocking Number]
Let $x,y\in M$ be two (not necessarily distinct) points in $(M,g)$.  The \textit{blocking number} $b_{g}(x,y)$ for $L_{g}(x,y)$ is defined by 
$$ b_{g}(x,y)=\inf \{n \in \mathbb{N} \cup \{\infty\} | \,L_{g}(x,y) \textup{ is blocked by $n$ points} \}. $$
\end{Def}

\begin{Def}[Finite Blocking Property]
A compact Riemannian manifold $(M,g)$ is said to have \textit{finite blocking} if $b_{g}(x,y)< \infty$ for every $(x,y)\in M \times M$.  When $(M,g)$ has finite blocking and the blocking numbers are uniformly bounded above, $(M,g)$ is said to have \textit{uniform} finite blocking. 
\end{Def}

The finite blocking property seems to have originated in the study of polygonal billiard systems and translational surfaces (see e.g. \cite{Fo}, \cite{Gu1}, \cite{Gu2}, \cite{Gu3}, \cite{HiSn}, \cite{Mo1}, \cite{Mo2}, \cite{Mo3}, and \cite{Mo4}).     Our motivation comes from the following theorem (see e.g. \cite{Fo}, \cite[lemma 1]{Gu1}, or \cite[proposition 2]{GuSc}):

\vskip 5pt

\noindent {\bf Theorem:}  Compact flat Riemannian manifolds have uniform finite blocking.

\vskip 5pt

We believe the following is true:

\begin{Conj}
Let $(M,g)$ be a compact Riemannian manifold with finite blocking.  Then $g$ is a flat metric.
\end{Conj}

There is a natural analogue of (uniform) finite blocking for general geodesic metric spaces.  We provide an extension of the above theorem in section 5:

\begin{Thm2}
Compact quotients of Euclidean buildings have uniform finite blocking.
\end{Thm2}

 As evidence for the above conjecture we prove the following theorem in section 4:

\begin{Thm2}
Let $(M,g)$ be a compact nonpositively curved Riemannian manifold with the finite blocking property.  Then $g$ is a flat metric.
\end{Thm2}

This theorem is a consequence of a well known result about nonpositively curved manifolds and the next theorem relating the finite blocking property to the topological entropy of the geodesic flow:

\begin{Thm2}
Let $(M,g)$ be a compact Riemannian manifold without conjugate points.  If $h_{top}(g)>0$, then $b_g(x,y)= \infty$ for every $(x,y) \in M$.  In other words, given any pair of points $x,y \in M$ and a finite set $F\subset M-\{x,y\}$, there exists a geodesic segment connecting $x$ to $y$ and avoiding $F$.
\end{Thm2}

Working independently and using a similar approach K. Burns and E. Gutkin have also obtained Theorem 3 as well as the following \cite{BuGu}:  

\vskip 5pt

\noindent {\bf Theorem:} (Burns and Gutkin)  Let $(M,g)$ be a compact Riemannian manifold with uniform finite blocking.  Then $h_{top}(g)=0$ and the fundamental group of $M$ is virtually nilpotent.  
\vskip 5pt

In section 2, we define \textit{regular finite blocking} by imposing a continuity and separation hypothesis on blocking sets.  We show that manifolds with regular finite blocking have uniform finite blocking and are conjugate point free.  Combining this result with the previous theorem of K. Burns and E. Gutkin, and recent work of N.D. Lebedeva \cite{Le} yields: 

\begin{Thm2}
Let $(M,g)$ be a compact Riemannian manifold with regular finite blocking.  Then $g$ is a flat metric.  
\end{Thm2}

Blocking light is also interesting in the context of the nonnegatively curved compact type locally symmetric spaces.  In \cite{GuSc}, they show the following:

\vskip 5pt

\noindent {\bf Theorem:} (Gutkin and Schroeder)
Let $(M,g)$ be a compact locally symmetric space of compact type with $\mathbb{R}$-rank $k\ge1$.  Then $b_g(x,y)\leq 2^k$ for almost all $(x,y) \in M \times M$.

\vskip 5pt

We refer the reader to \cite{GuSc} for a more precise formulation and discussion of this result.  On the positively curved side, we focus on the blocking properties of the compact rank one symmetric spaces or CROSSes.  The CROSSes are classified and consist of the round spheres $(S^n, can)$, the projective spaces $(K\mathbb{P}^n,can)$ where $K$ denotes one of $\mathbb{R}$,$\mathbb{C}$, or $\mathbb{H}$, and the Cayley projective plane $(Ca\mathbb{P}^2, can)$ where $can$ denotes a symmetric metric.  The CROSSes all satisfy the following blocking property:

\begin{Def}[Cross Blocking]
A compact Riemannian manifold $(M,g)$ is said to have cross blocking if $$0<d(x,y)<\Diam(M,g) \implies b_g(x,y)\leq 2.$$
\end{Def}

Just as with finite blocking, we also define \textit{regular cross blocking} by imposing a continuity and separation hypothesis on blocking sets.  In addition to cross blocking, round spheres also satisfy the following blocking property:

\begin{Def}[Sphere Blocking]
A compact Riemannian manifold $(M,g)$ is said to have sphere blocking if $b_g(x,x)=1$ for every $x \in M$.
\end{Def}

This is a blocking interpretation of  ``antipodal points''; we think of the single blocker for $L_g(x,x)$ 
as being antipodal to $x$.  We believe the following is true:

\begin{Conj}
A compact Riemannian manifold $(M,g)$ has cross blocking if and only if $(M,g)$ is isometric to a compact rank one symmetric space.  In particular, $(M,g)$ has cross blocking and sphere blocking if and only if $(M,g)$ is isometric to a round sphere.
\end{Conj}

As support for this conjecture,  we prove the following theorems in section 3:

\begin{Thm2}
Let $(S^2,g)$ be a metric on the two sphere with cross blocking and sphere blocking.  Then a shortest periodic geodesic is simple with period $2\Diam(S^2,g)$.
\end{Thm2}


\begin{Thm2}
Let $(M^{2n},g)$ be an even dimensional manifold with positive sectional curvatures and regular cross blocking.  Then $(M,g)$ is a Blaschke manifold.  If in addition $M$ is diffeomorphic to a sphere or has sphere blocking, then $(M,g)$ is isometric to a round sphere.
\end{Thm2}

\begin{Thm2}
Let $(M,g)$ be a compact Riemannian manifold with regular cross blocking, sphere blocking, and which doesn't admit a nonvanishing line field.  Then $(M,g)$ is isometric to an even dimensional round sphere. 
\end{Thm2}

\vskip 5pt

\centerline {\bf Acknowledgements:}  

\vskip 5pt

We owe our gratitude to Keith Burns; Keith pointed out a subtlety regarding Theorem 3 (Proposition 4.2 below) that we originally missed during the earlier stages of this work.  We also thank Ralf Spatzier for numerous helpful discussions.  The first author would like to thank the \'Ecole Normale Sup\'erieure (Lyon) for their hospitality.  The work of the first author was partially supported by the National Science Foundation under grant DMS-0606002.  The second author would like to thank The Ohio State University mathematics department for their hospitality.  This research was partially conducted during the period the second author was employed by the Clay Mathematics Institute as a Liftoff Fellow.

\section{Finite Blocking and Conjugate Points}

For the reader's convenience, we begin with the definition of a conjugate point in a compact Riemannian manifold $(M,g)$.  We let $TM$ (resp. $UM$) denote the tangent bundle (resp. unit tangent bundle) of $M$ and denote the fibers above a point $p\in M$, by $T_{p}M$ and $U_{p}M$.  For a point $p\in M$, the exponential map $$\exp_{p}:T_{p}M \rightarrow M$$ is everywhere defined by completeness.  

\begin{Def}[Conjugate Point]
A point $q=\exp_{p}(v) \in M$ is \textit{conjugate} to $p$ along the unit speed geodesic $\gamma_{v}: [0, ||v||] \rightarrow M$ with initial condition $\frac{v}{||v||} \in U_{p}(M) \subset UM$ if $d(\exp_{p})_{v}$ is not of full rank.
\end{Def}

In \cite{Wa1}, F. Warner describes the conjugate locus of singular points $C(p)\subset T_p(M)$ for the exponential map $\exp_p$.  A point $v\in C(p)$ is said to be \textit{regular} if there exists a neighborhood $U$ of $v$ such that each ray emanating from the origin in $T_p(M)$ intersects at most one point in $C(p) \cap U$.  The \textit{order} of a point $v \in C(p)$ is defined to be the dimension of the kernel of $d(\exp_p)_v$.  Warner shows that the set of regular points $C^R(p) \subset C(p)$ is an open dense subset of $C(p)$ which (if nonempty) forms a codimension one submanifold of $T_p(M)$.  Moreover, the order of points is constant in each connected component of $C^R(p)$ and there are normal forms depending on the order of the point for the exponential map in a neighborhood of each regular point.  From these normal forms, it follows that the preimage under the exponential map of a regular conjugate point of order more than one is indiscrete.  It appears that regular points of order more than one are rare in Riemannian manifolds (see e.g. \cite{Wa2}).  The next proposition shows that there are no such conjugate points in Riemannian manifolds with the finite blocking property.

\begin{Prop}
Suppose that $(M,g)$ is a compact Riemannian manifold with finite blocking.  Then for each $p\in M$, point preimages of $\exp_p$ are discrete subsets of $T_pM$.
\end{Prop}

\begin{proof}
Suppose not. Then there are points $p,q \in M$ and a sequence of vectors $\{v_i\} \subset \exp_p^{-1}(q)$  converging to a vector $v_{\infty} \in \exp_p^{-1}(q)\subset T_p(M)$.  Let $B=\{b_1,\ldots,b_k\} \subset M$ be a finite blocking set for $L_{g}(p,q)$.  Define first blocking times $t_i \in (0,1)$ by $$t_i=\inf\{t\in (0,1)| \exp_p(tv_i) \in B \}.$$  After possibly relabeling blockers and passing to a subsequence, we may assume that $\exp_p(t_i v_i)=b_1$ for all $i\in \mathbb{N}$.  A subsequence of the vectors $\{t_i v_i\}$ converge to a vector $t_{\infty} v_{\infty}$ and by continuity of the exponential map, $\exp_p(t_{\infty} v_{\infty})=b_1$.  This shows that the point $b_1$ is a sooner conjugate point to $p$ along the geodesic ray $\gamma(t)=\exp_p(tv_{\infty})$ than is the point $q$.  By repeating this argument, there is always a sooner conjugate point, contradicting the fact that conjugate points are discrete along a geodesic.
\end{proof}

We expect that the light between conjugate points will never be finitely blocked.  Next we impose some restrictions on blocking sets and show that the light between conjugate points cannot be finitely blocked by such sets.  For a compact Riemannian manifold $(M,g)$, let $M' \subset M \times M$ be the subset of points for which $b_g < \infty$, $T'\subset TM$ be the subset of vectors $(p,v)\in TM$ for which $(p,\exp_p(v))\in M'$, and let $\mathcal{F}(M)$ denote the set of finite subsets of $M$.  A \textit{blocking function} for $(M,g)$ is a symmetric map $\mathcal{B}: M' \rightarrow \mathcal{F}(M)$ such that for each $(x,y) \in M',$ $\mathcal{B}(x,y)$ is a finite blocking set for $L_{g}(x,y)$.  Given a blocking function $\mathcal{B}$ we define the \textit{first blocking time} $t_{\mathcal{B}}:T' \rightarrow (0,1)$ by $$t_{\mathcal{B}}(p,v)=\inf\{t \in (0,1)| \exp_p(tv) \in \mathcal{B}(p,\exp_p(v))\}.$$

\begin{Def}[Continuous Blocking]
We say that a closed Riemannian manifold $(M,g)$ has \textit{continuous blocking} if there is a blocking function $\mathcal{B}$ for which the first blocking time $t_{\mathcal{B}}:T' \rightarrow (0,1)$ is continuous.\end{Def}

\begin{Def}[Separated Blocking]
We say that a blocking function $\mathcal{B}$ is \textit{separated} if there exists an $\epsilon>0$ such that the $\epsilon$-neighborhoods of blocking points in each finite blocking set $\mathcal{B}(x,y)\subset M$ are disjoint.
\end{Def}


\begin{Lem}
Let $\mathcal{B}$ be a separated blocking function for a compact Riemannian manifold $(M,g)$.  Then the cardinalities of blocking sets defined by $\mathcal{B}$ are uniformly bounded above.  In particular, compact Riemannian manifolds with finite blocking and a separated blocking function have uniform finite blocking.   
\end{Lem}

\begin{proof}
Suppose that the blocking sets defined by $\mathcal{B}$ are $\epsilon$-separated.  An upper bound $K_{max}$ for the sectional curvatures yields a lower bound $C:=C(K_{max},\epsilon)>0$ for the volume of balls of radius $\epsilon$ in $M$.  Therefore, there are at most $\vol(M,g)/C$ disjoint balls of radius $\epsilon$ in $M$, concluding the proof.
\end{proof}

\begin{Prop}
Let $(M,g)$ be a compact Riemannian manifold with a blocking function $\mathcal{B}$ that is both continuous and separated.  Let $p \in M$ and suppose that $U$ is an open subset of  $T_p(M)$ consisting of vectors from $T' \cap T_p(M)$.  Then $p$ is not conjugate to any point in $\exp_p(U) \subset M$.    
\end{Prop} 

\begin{proof}
Suppose not.  Then there is a vector $v \in U$ for which $\exp_p(v)$ is conjugate to $p$.  It is well known (see e.g. \cite{Wa1}) that $\exp_p$ is not one to one in any neighborhood of $v$.  Let $B_i$ be a sequence of balls centered at $v$ and contained in $U$ with radii decreasing to zero.  For each $i$, choose distinct points $x_i,y_i \in B_i$ with $\exp_p(x_i)=\exp_p(y_i):=q_i$.  Let  $\mathcal{B}$ be a continuous and separated blocking function, and let $l_i:=\exp_p(t_{\mathcal{B}}(x_i)x_i)$ and $r_i:=\exp_p(t_{\mathcal{B}}(y_i)y_i)$ be the associated first blocking points in $\mathcal{B}(p,q_i)$.  By continuity of $\mathcal{B}$, the sequences $\{l_i\}$ and $\{r_i\}$ both converge to $\exp_p(t_{\mathcal{B}}(v)v)$.  This contradicts the separatedness of $\mathcal{B}$ for all large enough indices $i \in \mathbb{N}$. 
\end{proof}

When a compact Riemannian manifold $(M,g)$ has finite (resp. cross) blocking and a continuous and separated blocking function, we shall say that $(M,g)$ has \textit{regular finite} (resp. \textit{regular cross}) \textit{blocking}.

\begin{Cor}
Let $(M^n,g)$ be a compact Riemannian manifold with regular finite blocking.  Then $(M,g)$ has uniform finite blocking and is conjugate point free.  In particular, the universal cover of $M$ is diffeomorphic to $\mathbb{R}^n$ 
\end{Cor}

\begin{proof}
Since $(M,g)$ has finite blocking, $T'=TM$.  By lemma 2.2, $(M,g)$ has uniform finite blocking and by proposition 2.3, $(M,g)$ is conjugate point free.  The second statement is Hadamard's theorem.
\end{proof}

\begin{Cor}
Let $(M,g)$ be a compact Riemannian manifold with regular cross blocking.  If $p,q \in M$ are conjugate points, then $d(p,q)=0$ or $d(p,q)=\Diam(M,g)$.
\end{Cor}

\begin{proof}
Suppose $p,q \in M$ satisfy $0<d(p,q)<\Diam (M,g)$ and let $\gamma:[0,1] \rightarrow M$ be a geodesic with $\gamma(0)=p$ and $\gamma(1)=q$.  Then $\gamma(t)=\exp_p(tv)$ for some $v \in T_p(M).$ By continuity of the exponential map and the cross blocking property, there is an open set $U\subset T_p(M)$ containing $v$ and satisfying $U \subset T' \cap T_p(M)$.  By proposition 2.3, the points $p$ and $q$ are not conjugate.
\end{proof}

\section{Blocking Light and Round Spheres} 
In this section, we show under various hypothesis that compact Riemannian manifolds with blocking properties similar to those of round spheres are necessarily isometric to round spheres.  In general, we believe the following should be true:

\begin{Conj}
Let $(M,g)$ be a compact Riemannian manifold with cross blocking.  Then $(M,g)$ is isometric to a compact rank one symmetric space.
\end{Conj}

We begin by reviewing the definition and basic properties concerning cut points in a compact Riemannian manifold $(M,g)$.  In this section, a unit speed geodesic $\gamma:[0,L_{\gamma}] \rightarrow M$ for which $\gamma(0)=\gamma(L_{\gamma})=p$ will be called a \textit{geodesic lasso} based at $p$.  For a geodesic lasso $\gamma$, we shall denote by $\gamma^{-1}:[0,L_{\gamma}] \rightarrow M$ the geodesic lasso obtained by traversing $\gamma$ in the reverse direction; specifically, $\gamma^{-1}(t):=\gamma(L_{\gamma}-t)$.  If in addition the geodesic $\gamma$ is regular at $p$, i.e. $\dot{\gamma}(0)=\dot{\gamma}(L_{\gamma})$, $\gamma$ will be called a \textit{closed geodesic} based at $p$.  By a \textit{simple lasso} (resp. \textit{simple closed geodesic})based at $p$ we mean a lasso (resp. closed geodesic) $\gamma:[0,L_{\gamma}] \rightarrow M$ based at $p$ which is injective on the interval $(0,L_{\gamma})$ and with $p \notin \gamma((0,L_{\gamma}))$.

\begin{Def}
Let $(M,g)$ be a compact Riemannian manifold, $p \in M$, $v\in U_p(M)$, and $\gamma:[0,\infty) \rightarrow M$ the unit speed geodesic ray defined by $\gamma(t)=\exp_p(tv)$ .  Let $[0,t_0]$ be the largest interval for which $t \in [0,t_0]$ implies $d(p,\gamma(t))=t$.  The point $\gamma(t_0)$ is said to be a cut point to $p$ along the geodesic $\gamma$.  The union of the cut points to $p$ along all the geodesics starting from $p$ is called the cut locus and will be denoted by $\Cut(p)$.   
\end{Def}

The next two propositions are well known and describe points in the cut locus (see e.g. \cite{do Ca}).

\begin{Prop}
Suppose that $\gamma(t_0)$ is the cut point of $p=\gamma(0)$ along a geodesic $\gamma$.  Then either: 

\begin{itemize}
\item $\gamma(t_0)$ is the first conjugate point of $\gamma(0)$ along $\gamma$, or
\item  there exists a geodesic $\sigma \neq \gamma$ from $p$ to $\gamma(t_0)$ such that $\length(\sigma)=\length(\gamma)$.
\end{itemize}

Conversely, if (a) or (b) is satisfied, then there exists  $t' \in (0,t_0]$ such that $\gamma(t')$ is the cut point of $p$ along $\gamma$.
\end{Prop}

\begin{Prop}
Let $p \in M$ and suppose that $q\in  \Cut(p)$ satisfies $d(p,q)=d(p,\Cut(p))$.  Then either:

\begin{itemize}
\item there exists a minimizing geodesic $\gamma$ from $p$ to $q$ along which $q$ is conjugate to $p$, or

\item there exist exactly two minimizing geodesics $\gamma$ and $\sigma$ from $p$ to $q$ that together form a simple geodesic lasso based at $p$ of length $2d(p,\Cut(p))$.
\end{itemize}

\end{Prop}

It follows from proposition 3.1 that $\exp_p$ is injective on a ball of radius $r$ centered at the origin in $T_p(M)$ if and only if $r<d(p,\Cut(p))$.  

\begin{Def}
The \textit{injectivity radius} of $(M,g)$ is defined to be $$\inj(M,g)=\inf_{p\in M} d(p, \Cut(p)).$$
\end{Def}

Note that the injectivity radius of a compact Riemannian manifold is never larger than its diameter.  Compact manifolds for which the injectivity radius equals the diameter are known as \textit{Blaschke manifolds}.  All of the compact rank one symmetric spaces are Blaschke and the well known Blaschke conjecture asserts that these are the only Blaschke manifolds.  We will use the following theorem from \cite{Ber1}, extending earlier work of \cite{Gr}:

\begin{Thm}[Berger]
Let $(M,g)$ be a Blashke metric on a smooth sphere.  Then the metric $g$ is a symmetric metric.
\end{Thm} 

For a Blaschke manifold $(M,g)$, all infinite geodesics $\gamma:\mathbb{R} \rightarrow M$ cover simple closed geodesics of length $2\Diam(M,g)$ (see e.g. \cite[Corollary 5.42]{Be}).  The next proposition is a first step in showing that geodesics in a manifold with cross blocking behave similarly to those in a Blaschke manifold. 

\begin{Prop}
Suppose that $(M,g)$ has cross blocking.  Let $\gamma:[0,L_{\gamma}] \rightarrow M$ be a unit speed simple geodesic lasso based at $p \in M$.  Then $L_{\gamma} \leq 2\Diam(M,g)$ and the point $\gamma(L_{\gamma}/2)$ is the cut point to $p$ in both of the directions $\dot \gamma(0)$ and $-\dot \gamma(L_{\gamma})$. 
 \end{Prop}

\begin{proof}
Let $D:=\Diam(M,g)$ and let $c_1$ be the cut point to $p$ in the direction $\dot \gamma(0)$ and $c_2$ be the cut point to $p$ in the direction $-\dot \gamma(L_{\gamma})$.  By simplicity of $\gamma$, there exists unique $t_1 , t_2 \in (0,L_{\gamma})$ such that $c_i=\gamma(t_i)$ for $i=1,2$.  By proposition 3.1, $t_1\leq L_{\gamma}/2 \leq t_2$.  The statement of the proposition follows from showing that $t_1=t_2$,  i.e that $c_1=c_2$ and the point $\gamma(L_{\gamma}/2)$ is the cut point to $p$ in both of the directions $\dot \gamma(0)$ and $-\dot \gamma(L_{\gamma})$.

We now will assume that $t_1<t_2$ and must obtain a contradiction.  With this assumption, we first argue that $d(p,\gamma((t_1,t_2)))=D$.  If not,  choose a $t \in (t_1,t_2)$ for which $0<d(p, \gamma(t))<D$.  Note that since $\gamma$ is simple, the restrictions of $\gamma$ to the interval $[0,t]$ and $\gamma^{-1}$ to the interval $[0,L_{\gamma}-t]$ define distinct elements in $L_{g}(p,\gamma(t))$ with nonintersecting interiors.  There must be a single blocking point on the restriction of $\gamma$ to the interior of each of these intervals since $(M,g)$ has cross blocking. As neither of these light rays are minimizing, there is a unit speed minimizing geodesic $\sigma:[0,L_{\sigma}] \rightarrow M$ joining $p$ to $\gamma(t)$.   Since $\sigma$ is minimizing, it defines a third light ray between $p$ and $\gamma(t)$, whence its interior must intersect one of the two blocking points for $L_{g}(p, \gamma(t))$.  Let $s':=\inf \{s\in (0,L_{\sigma})| \, \sigma(s) \in \interior{\gamma}!
  \}$.  By simplicity of $\gamma$, there is a unique $t' \in (0,L_{\gamma})$ such that $\sigma(s')=\gamma(t'):=q$. Since $\sigma$ is unit speed and minimizing, $0<d(p,q)=s'<L_{\sigma}=d(p,\gamma(t))<D$, so that $b_g(p,q)\leq 2$ by the cross blocking condition.  However, the restriction of $\gamma$ to $[0,t']$, $\gamma^{-1}$ to $[0,L_{\gamma}-t']$, and the restriction of $\sigma$ to $[0,s']$ define three distinct elements in $L_{ g}(p,q)$ with nonintersecting interiors, implying $b_g(p,q) \geq 3$. This is a contradiction,  whence $d(p, \gamma((t_1,t_2)))=D$.

Next we show that $d(p,\gamma((t_1,t_2)))=D$ yields a contradiction, completing the proof.  By the discreteness of conjugate points along geodesics, there is a $t \in (t_1,t_2)$ so that $q:=\gamma(t)$ is not conjugate to $p$ in the direction $\dot \gamma(0)$ or in the direction $-\dot \gamma(L_{\gamma})$.  Neither of the restrictions of $\gamma$ to $[0,t]$ or $\gamma^{-1}$ to $[0,L_{\gamma}-t]$ is minimizing so that there is a unit speed minimizing geodesic $\sigma:[0,D]\rightarrow M$ joining $p$ to $q$.  Note that since $\sigma$ is minimizing, the interior of $\sigma$ cannot intersect $\gamma$. Indeed, a first point of intersection between the interiors of $\sigma$ and $\gamma$ would be at a point $p'$ satisfying $0<d(p,p')<D$ so that the reasoning from the previous paragraph may be applied to obtain a contradiction.  Let $q_n=\sigma(D-1/2^n)$.  Choose sufficiently small neighborhoods $B_1$ of $t\dot \gamma(0)$ and $B_2$ of 
$-(L_{\gamma}-t)\dot \gamma(L_{\gamma})$ on which $\exp_p$ restricts to a local diffeomorphism.  For all sufficiently large $n$, there are unique $x_n \in B_1$ and $y_n \in B_2$ such that $q_n=\exp_p(x_n)=\exp_p(y_n)$.  It follows by the continuity properties of the exponential map that for suitably large $n$, the geodesics $$s \mapsto \exp_p(sx_n)$$ and $$s \mapsto \exp_p(sy_n)$$ for $s\in[0,1]$ and the restriction of $\sigma$ to the interval $[0,D-1/2^n]$ define three light rays between $p$ and $q_n$ with nonintersecting interior.  Hence, $b_g(p,q_n) \geq 3$ for suitably large $n$.  But $0<d(p,q_n)<D$, so that $b_g(p,q_n) \leq 2$ by cross blocking, a contradiction.

\end{proof}

\begin{Lem}
Suppose that $(M,g)$ has cross blocking and sphere blocking. Suppose that $\gamma:[0,L_{\gamma}] \rightarrow M$ is a unit speed simple geodesic lasso based at $p\in M$.  If $L_{\gamma}<2\Diam(M,g)$, then $\gamma$ is regular at $p$ and all lassos based at $p$ finitely cover $\gamma$.  If $L_{\gamma}=2\Diam(M,g)$, then the interior of all of the geodesic lassos through $p$ intersect in the point $\gamma(L_{\gamma}/2)$. 
\end{Lem} 

\begin{proof}
Let $\overline{p}:=\gamma(L_{\gamma}/2)$.  Suppose there is a (not necessarily simple) unit speed lasso $\tau:[0,L_{\tau}]\rightarrow M$ through $p$ with $\dot \tau(0)$ distinct from $\dot \gamma(0)$ and $-\dot \gamma (L_{\gamma})$.  As $(M,g)$ has sphere blocking and $\gamma$ is simple, there is a unique $t\in (0,L_{\gamma})$ such that $\gamma(t)$ blocks $L_{g}(p,p)$.  Let $s:=\inf\{t \in (0,L_{\tau}]\, | \, \tau(t)=p\}$. The restriction of $\tau$ to the interval $[0,s]$ gives an element in $L_g(p,p)$ so that by sphere blocking, its interior must pass through the blocking point $\gamma(t)$ (and hence $\interior(\gamma)$).  Let $s':=\inf \{t \in (0,s)|\,  \tau(t) \in \interior{\gamma} \}$. By simplicity of $\gamma$ there is a unique $t' \in (0,L_{\gamma})$ such that $\gamma(t')=\tau(s'):=q$.  As $\gamma$ is simple, the restrictions of $\gamma$ to the intervals $[0,t']$, $\gamma^{-1}$ to the interval $[t',L_{\gamma}]$, and 
$\tau$ to the interval $[0,s']$ define three distinct light rays between $p$ and $q$ with nonintersecting interiors. Since $p\neq q$ cross blocking implies that $d(p,q)=\Diam(M,g)$ and that $q=\overline{p}$  (by proposition 3.4 and since $L_{\gamma}/2 \leq \Diam(M,g)$).  Hence, if $L_{\gamma}/2<\Diam(M,g)$ there are no geodesic lassos through $p$ with initial tangent vector outside of the set $\{\dot\gamma(0), -\dot \gamma(L_{\gamma}) \}$, and if $L_{\gamma}/2=\Diam(M,g)$, the interior of any lasso through the point $p$ passes through the point $\overline{p}$.  This concludes the proof of the last statement in the lemma. 

We now assume that $L_{\gamma}/2< \Diam(M,g)$, and will argue that $\gamma$ is regular at the point $p$ and that all lassos at $p$ finitely cover $\gamma$.  By simplicity of $\gamma$, the restriction of $\gamma$ to the intervals $[0,L_{\gamma}/2]$ and $\gamma^{-1}$ to the interval $[0,L_{\gamma}/2]$ define distinct light rays between $p$ and $\overline{p}$ with nonintersecting interiors.  As $(M,g)$ has cross blocking, the interior of a third light ray must intersect the interior of $\gamma$ in a blocker and will therefore have a first point of intersection $p'$ with the interior of $\gamma$.  This implies $b_g(p,p')\geq 3$, a contradiciton.   Therefore $|L_g(p,\overline{p})|=2$, while $G_g(p,\overline{p})$ is infinite by \cite{Se}.  Note that any geodesic segment from  $G_g(p,\overline{p})-L_g(p,\overline{p})$ is obtained from extending one of the two light rays in $L_g(p, \overline{p})$.  Each such extension  gives rise to a geodesic lasso based at $p$ with initial tangent vector in the set $\{-\dot\gamma(0),\dot\gamma(L_{\gamma})\}$.  But by the previous paragraph, the initial tangent vector of all lassos at $p$ lie in the set $\{\dot\gamma(0),-\dot\gamma(L_{\gamma})\}.$  Therefore, $\gamma$ must be regular at $p$ and all lassos at $p$ finitely cover $\gamma$.  
\end{proof}

\begin{Def}
A $SC_{2L}$ manifold is a Riemannian manifold with the property that all geodesics cover simple closed geodesics of length $2L$.  
\end{Def}

It is tempting to think that the only $SC_{2L}$ manifolds are the CROSSes.  Amazingly, O. Zoll exhibited an exotic $SC_{2L}$ real analytic Riemannian metric on the two sphere \cite{Zo}.  This example is discussed in \cite[Chapter 4]{Be} along with examples on higher dimensional spheres.  We remark that these examples are not cross blocked.  Indeed, proposition 3.4 implies that $SC_{2L}$ manifolds with cross blocking are Blaschke with $\inj(M,g)=\Diam(M,g)=L$, while theorem 3.3 asserts that there are no exotic Blaschke metrics on spheres.  
In view of our conjecture that the manifolds with cross blocking are precisely the CROSSes and the Blaschke conjecture that the Blaschke manifolds are precisely the CROSSes, we expect that Blaschke manifolds are precisely those manifolds with cross blocking.  In the next proposition, we use well known results concerning Blaschke manifolds to show that Blaschke manifolds all have cross blocking.

\begin{Prop}
Suppose that $(M,g)$ is a Blaschke manifold.  Then $(M,g)$ has cross blocking.
\end{Prop}

\begin{proof}
Suppose that $p,q \in M$ satisfy $0<d(p,q)<D:=\Diam(M,g)$.  By \cite[corollary 5.42]{Be}, $(M,g)$ is a $SC_{2D}$ manifold.   It follows that there is a unit speed simple closed geodesic $\gamma:[0,2D] \rightarrow M$ with $\gamma(0)=p$ and $\gamma(d(p,q))=q$.  The restriction of $\gamma$ to the intervals $[0,d(p,q)]$ and $\gamma^{-1}$ to $[0,2D-d(p,q)]$ give two distinct elements in $L_{g}(p,q)$ with nonintersecting interiors.  Hence $b_g(p,q) \geq 2$.  If $b_g(p,q)>2$ then there must be a third unit speed light ray $\beta:[0,L_{\beta}] \rightarrow M$ joining $p$ to $q$.  Note that since $\beta$ is a light ray and since all geodesics are periodic with period $2D$, $L_{\beta}<2D$.  Moreover, since $d(p,q)<D=\inj(M,g)$, the restriction of $\gamma$ to the interval $[0,d(p,q)]$ is the unique minimizing geodesic from $p$ to $q$ so that $L_{\beta}>\inj(M,g)=D$.  Extend $\beta$ to the simple closed geodesic $\overline{\beta}:[0,2D] \rightarrow M$.  Then the restriction of $\overline{\beta}$ to the interval $[L_{\beta},2D]$ gives a geodesic joining $q$ to $p$ of length $2D-L_{\beta}<D=\inj(M,g)$.  Hence, there are two minimizing geodesics joining $p$ and $q$, a contradiction.  Therefore, $b_g(p,q)=2$ and Blaschke manifolds have cross blocking.  
\end{proof}

\begin{Cor}
Suppose that $(M,g)$ is Blaschke manifold with sphere blocking.  Then $(M,g)$ is isometric to a round sphere.
\end{Cor}

\begin{proof}
By the last proposition $(M,g)$ has cross blocking.  By lemma 3.5, the interior of all of the simple closed geodesics through $p$ intersect in a single point $p'$ satisfying $d(p,p')=\Diam(M,g)$.  By proposition 3.4, $p'$ is the cut point to $p$ along all these geodesics.  Hence, $\dim(\Cut(p))=0$, from which it follows (see e.g. \cite[proposition 5.57]{Be}) that $M$ is diffeomorphic to a sphere.  By theorem 3.3, $(M,g)$ is a round sphere. 
\end{proof}

For the proof of the next theorem, we will need the following two definitions:

\begin{Def}
For $p\in M$ and $p' \in \Cut(p)$, define the link $\Lambda(p,p') \subset U_{p'}M$ by  $$\Lambda(p,p')=\{-\dot{\gamma}(d(p,p'))|\,  \gamma \, \text{is a unit speed and minimizing geodesic from $p$ to $p'$}\}.$$
\end{Def}

\begin{Def}
For $p \in M$ and $U \subset U_{p}M$, $U$ is said to be a great sphere if $U$ is the intersection of a linear subspace of $T_{p}M$ with $U_{p}M$.
\end{Def}

\begin{Thm}
Suppose that $(M,g)$ is a compact Riemannian manifold with regular cross blocking, sphere blocking, and which does not admit a nonvanishing continuous line field. Then $(M,g)$ is isometric to an even dimensional round sphere.
\end{Thm}

\begin{proof}
First we argue that $(M,g)$ is a Blaschke manifold.  

To obtain a contradiction, suppose that $\inj(M,g)<\Diam(M,g):=D$.  We begin by showing that for each point $x \in M$ satisfying $d(x,\Cut(x))<D$ there is a unique simple closed geodesic based at $x$ and this geodesic has length $2d(x,\Cut(x))$.  Indeed, let $x$ satisfy $d(x,\Cut(x))<D$ and $x' \in \Cut(x)$ satisfy $d(x,x')=d(x,\Cut(x))$.  By corollary 2.5 the points $x$ and $x'$ are not conjugate so that by proposition 3.2 there is a simple geodesic lasso $C$ through $x$ of length $2\inj(M,g)$.  By lemma 3.5, $C$ is a simple closed geodesic through $x$ and is the unique lasso through $x$, as required.

Let $L:=\sup\{d(p,\Cut(p))\,| \, p\in M\}\leq D$.  If $L<D$, the preceeding paragraph shows that there is a unique closed geodesic $C_p$ of length less than $2D$ through each point $p \in M.$  Since all of the $C_q$ have lengths uniformly bounded above, whenever a sequence of points $\{p_i\}$ converge to a point $p_{\infty} \in M$, the sequence of closed geodesics $C_{p_{i}}$ converge to a closed geodesic $C_{\infty}$. By the uniqueness of these geodesics, $C_{\infty}=C_{p_{\infty}}$.  Therefore, the tangent spaces to these geodesics define a nonvanishing continuous line field on $M$, a contradiction.   

Hence, there is a point $p\in M$ satisfying $d(p,\Cut(p))=D$.  Such a point is said to have \textit{spherical cut locus} at $p$ \cite[definition 5.22]{Be}. By \cite[proposition 5.44]{Be}, the link $\Lambda(p,p') \subset U_{p'}(M)$ is a great sphere for each $p' \in \Cut(p)$, whence all geodesics through $p$ are periodic of period $2D$.  Now consider a geodesic connecting $p$ to a point $q$ satisfying $d(q,\Cut(q))<D$.  This geodesic gives rise to a closed geodesic of length $2D$ through $q$, while the first paragraph shows that there is a closed geodesic through $p$ of length $2d(q,\Cut(q))$.  This contradicts lemma 3.5, 
implying that at every point $q\in M$, we have $d(q,\Cut(q))=D$, and hence concluding the proof that  
$(M,g)$ is Blaschke.  

By corollary 3.7, $(M,g)$ is isometric to a round sphere.  As $M$ does not admit a nonvanishing line field, $(M,g)$ is isometric to an even dimensional round sphere.  
\end{proof}

Next, we adapt Klingenberg's estimate on the injectivity radius to obtain the following (see e.g. \cite[chapter 13, proposition 3.4]{do Ca}):

\begin{Thm}
Suppose that $(M^{2n},g)$ is an even dimensional, orientable, Riemannian manifold with positive sectional curvatures.  If $(M,g)$ has regular cross blocking, then $(M,g)$ is Blaschke.  In particular, if $M$ is diffeomorphic to a sphere or $(M,g)$ has sphere blocking, then $g$ is a round metric on a sphere. 
\end{Thm}

\begin{proof}
Suppose to the contrary that $\inj(M,g)<\Diam(M,g)$ and choose $p,q \in M$ so that $q \in \Cut(p)$ and $d(p,q)=\inj(M,g)$.  By corollary 2.5, $p$ and $q$ are not conjugate points.  By proposition 3.2, there is a unit speed simple closed geodesic $C:[0,2\inj(M,g)] \rightarrow M$ of length $2\inj(M,g)$ passing through $p=C(0)$ and $q=C(\inj(M,g))$.  Since $M$ is orientable and even dimensional, parallel transport along $C$ leaves invariant a vector  $v$ orthogonal to $C$ at $C(0)$.  The field $v(t)$ along $C(t)$ is the variational field of closed curves $C_s(t)$ for $s\in [0,\epsilon)$.  As the sectional curvatures are strictly positive, the second variational formula implies that $\length(C_s)<\length(C)$ for all small $s>0$.  For each $s>0$, let $q_s$ be a point of $C_s$ at maximum distance from $C_s(0)$.  Necessarily, $\lim_{s\rightarrow 0} q_s =q$ and $d(q_s,C_s(0))<\inj(M,g)$.  For each $s>0$, let $\gamma_s$ be the unique minimizing geodesic joining $q_s$ to $C_s(0)$.  No!
 te that each 
$\dot{\gamma_s}(0)$ is orthogonal to $C_s$ by the first variational formula.  Let $w\in T_{q}M$ be an accumulation point of the vectors $\dot{\gamma_s}(0) \in T_{q_s}M$ and $\gamma:[0,1]\rightarrow M$ be the geodesic defined by $\gamma(t):=\exp_{q}(tw)$.  It follows that $\gamma$ is a minimizing geodesic joining $q$ to $p$ which is orthogonal to $C$ at $q$.  As $\gamma$ is minimizing, it cannot intersect $C$ except at the points $p$ and $q$, whence $b_g(p,q) \geq 3$, a contradiction.  Therefore $(M,g)$ is Blaschke.  The last statement follows from theorem 3.3 and corollary 3.7. 
\end{proof}

\begin{Thm}
Suppose that $(S^2,g)$ is a Riemannian metric on the two sphere with cross blocking and sphere blocking.  Then a shortest nontrivial closed geodesic is simple and has length $2\Diam(S^2,g)$.
\end{Thm}

\begin{proof}
Let $D:=\Diam(M,g)$ and let $C$ be a shortest nontrivial closed geodesic.  We first argue that if $C$ is simple, then its length is $2\Diam(S^2,g)$.  Indeed, by proposition 3.4, $\length(C)\leq 2D$.  We suppose that $\length(C)<2D$ and will obtain a contradiciton.  Note that $C$ separates $S^2$ into two components.  By Santalo's formula (see e.g. \cite{Sa} pg. 488 or \cite{Ber2} pg. 290), almost all geodesic rays with initial point on $C$ eventually leave the component they initially enter.  Choose one such ray $\gamma:[0,\infty) \rightarrow M$ and let $t:=\inf \{t \in (0,\infty)|\, \gamma(t)\in C\}.$  If $\gamma(t)$ is distinct from $\gamma(0)$, then $b_g(\gamma(0),\gamma(t))\geq 3$, a contradiction.  Hence, $\gamma(0)=\gamma(t)$, also a contradiction by lemma 3.5. 

Next we argue that a shortest nontrivial closed geodesic $C$ is simple.  By \cite{NaRo} or \cite{Sab}, $\length(C) \leq 4D$.  Suppose that $C$ is not simple, and choose a unit speed paramaterization $C:[0,L_C] \rightarrow M$ such that $C(0)$ is a crossing point.  Let $s=\inf \{ t \in (0,L_C) \, | \, C(s)=C(0)\}$.  Without loss of generality, we may assume that the restriction of $C$ to the interval $[0,s]$ defines a simple lasso at $C(0)$, whence $s=2D$ by proposition 3.4.  The restriction of $C$ to the interval $[2D,L_C]$ defines another lasso at $C(0)$.  If this lasso is not simple, then it contains a simple lasso of length $2D$ and $L_C>4D$, a contradiction.  Hence, the restriction of $C$ to the interval $[2D, L_C]$ defines a simple lasso and $L_C=4D$.  By lemma 3.5, $C(D)=C(3D)$.  Note that the restriction of $C$ to $[0,2D]$ separates $S^2$ into two components. This implies that $C((0,2D))\cap C((3D,4D)) \neq \emptyset$.  Letting $s=\inf \{t \in (3D,4D) \, | \, C(t) \in C((0,2D))\},$ it follows that $b_g(C(D),C(s)) \geq 3$, a contradiction. 
\end{proof}

\section{Finite Blocking Property and Entropy}

In this section we relate the finite blocking property for a compact Riemannian manifold $(M,g)$ to the topological entropy $h_{top}(g)$ of its geodesic flow.

Our starting point is a well known theorem (see e.g. \cite[corollary 1.2]{Ma})  identifying the topological entropy with the exponential growth rate of the number of geodesics between pairs of points in $M$.  For $x,y \in M$ and $T>0$, let $n_{T}(x,y)$ (resp. $m_{T}(x,y)$) denote the number of geodesic segments (resp. light rays) between the points $x$ and $y$ of length no more than $T$.

\begin{Thm}[Ma\~n\'e]
Let $(M,g)$ be a compact Riemannian manifold without conjugate points.  Then  $$h_{top}(g)=\lim_{T \rightarrow \infty} \frac{\log(n_{T}(x,y))}{T},$$ for all $(x,y) \in M \times M.$ 
\end{Thm}

The main observation of this section lies in the following:

\vskip 5pt

\begin{Prop}
Let $(M,g)$ be a compact Riemannian manifold without conjugate points.  If $h_{top}(g)>0$, then $b_{g}(x,y)= \infty$ for all $(x,y) \in M \times M.$
\end{Prop}

\begin{proof}

Let $I:=\inj(M,g)$.  We first argue that $n_T(x,y) \leq (T/2I)^2 \,m_T(x,y)$.  For a unit speed geodesic $\gamma:[0,L_{\gamma}] \rightarrow M$ in $G_{g}(x,y)$, let $t_1(\gamma):= \sup \{t \in [0,L_{\gamma}) \,|\, \gamma(t)=x\}$ and $t_2(\gamma):=\inf \{t \in(t_1(\gamma),L_{\gamma}] \, | \, \gamma(t)=y\}$.  Restrict $\gamma$ to the interval $[t_1(\gamma),t_{2}(\gamma)]$ and change the parameter of this interval to define a unit speed geodesic $$\Light(\gamma):[0, t_2(\gamma)-t_1(\gamma)]\rightarrow M$$ in $L_{g}(x,y)$.  Since for each $\beta \in L_{g}(x,y)$, $\Light(\beta)=\beta$, it follows that $$\Light:G_g(x,y) \rightarrow L_{g}(x,y)$$ defines a surjective map.  To conclude this step, it suffices to show that given a fixed unit speed $\beta \in L_{g}(x,y)$ of length not more than $T>0$, there are at most $(T/2I)^2$ distinct preimages of $\beta$ under the map $\Light$ of length not more than $T$.  For each unit speed geodesic $\gamma:[0,L_{\gamma}] \rightarrow M$ of length not more than $T$ satisfying $\Light(\gamma)=\beta$, $\dot \gamma (t_1({\gamma}))= \dot \beta(0)$ and $\dot \gamma( t_2(\gamma))=\dot \beta(L_{\beta})$ since geodesics are determined by their initial conditions.  It follows that the number of preimages of $\beta$ having length bounded above by $T$ coincides with the number of different extensions of $\beta$ to a unit speed geodesic $\overline{\beta} \in G_g(x,y)$ having length not more that $T$.  Given such an extension $\overline{\beta}$, let $n_{\overline{\beta}}(x)$ and $n_{\overline{\beta}}(y)$ be the number of returns to $x$ and the number or returns to $y$.  Necessarily, $n_{\overline{\beta}}(x),n_{\overline{\beta}}(y) \leq T/2I$ since each return to $x$ or return to $y$ increases the length of $\beta$ by at least $2I$.  Hence, by uniqueness of geodesics, $n_T(x,y) \leq (T/2I)^2 \,m_T(x,y)$.

To complete the proof, we argue by contradiction, assuming there is a pair of points $x,y \in M$ with a finite blocking set $F=\{b_1, \ldots b_{k}\}\subset M-\{x,y\}$ for $L_g(x,y)$.    By definition, the interior of any light ray $\gamma\in L_g(x,y)$ passes through some point $b_{i} \in F$, breaking $\gamma$ into two geodesic segments $\gamma_1 \in G_g(x,b_i)$ and $\gamma_2 \in G_g(b_i,y)$.  If $\gamma$ has length bounded above by $T$, then one of $\gamma_1$ or $\gamma_2$ must have length bounded above by $T/2$.  Moreover, given a geodesic segment $\alpha \in G_g(x,b_i)$ (resp. $\beta \in G_g(y,b_i)$), there is at most one extension of $\alpha$ (resp. $\beta$) to a light ray between $x$ and $y$.  It follows that $m_T(x,y) \leq \sum_{j=1}^{k} n_{T/2}(x,b_j) + n_{T/2}(b_j,y)$.  Combining this with the estimate from the previous paragraph yields:
$$n_T(x,y) \leq (T/2I)^2\, \sum_{j=1}^{k} n_{T/2}(x,b_j) + n_{T/2}(b_j,y).$$

Let $0< \epsilon < h_{top}(g)/3$.  By theorem 4.1, there is a $T_0 \in \mathbb{R}$ so that $T>T_{o}$ implies $$|h_{top}(g)-\frac{\log(n_{T}(*_1,*_2))}{T}|< \epsilon,$$ for all $*_1, *_2 \in \{x,y \}\cup F.$  Therefore  

$$\exp^{(h_{top}(g) - \epsilon)T}< n_T(*_1,*_2)<\exp^{(h_{top}(g)+\epsilon)T},$$ for all $T>T_{0}$ and $*_1,*_2 \in \{x,y\} \cup F$.  It now follows that $$\exp^{(h_{top}(g)-\epsilon)T} < n_T(x,y) < 2k(T/2I)^2 \exp^{(h_{top}(g)+ \epsilon)T/2},$$ a contradiction for all sufficiently large values of $T$.  
\end{proof}

We remark that the conclusion  $b_g(x,y)= \infty$ for all $(x,y) \in M \times M$ may be phrased more geometrically as saying that given any point $(x,y) \in M\times M$ and any finite set $F \subset M-\{x,y\}$, there is a geodesic segment between $x$ and $y$ avoiding $F$. 
\vskip 5pt

As a corollary of proposition 4.2, we obtain the following:

\begin{Thm}
Let $(M,g)$ be a compact Riemannian manifold with nonpositive sectional curvatures.  Then $(M,g)$ has finite blocking if and only if $(M,g)$ is flat.
\end{Thm}

\begin{proof}
Assume that $(M,g)$ has nonpositive curvature and is not flat. Then $(M,g)$ has no conjugate or focal points. By \cite[corollary 3]{Pe}, a geodesic flow on a nonflat compact Riemannian manifold without focal points has positive entropy.  By proposition 4.2, $(M,g)$ does not have finite blocking. 
\end{proof}

K. Burns and E. Gutkin \cite{BuGu} made the nice observation that by assuming uniform finite blocking and by iterating the line of reasoning used in the proof of proposition 4.2 one can establish the following:  
 
 \begin{Thm}[Burns-Gutkin]
 Let $(M,g)$ be a compact Riemannian manifold with the uniform finite blocking.  Then $h_{top}(g)=0$ and $\pi_1(M)$ has polynomial growth.
 \end{Thm}
 
 Using their result we obtain the following:
 
 \begin{Thm}
 Let $(M,g)$ be a compact Riemannian manifold with regular finite blocking.  Then $(M,g)$ is flat.
 \end{Thm}
 
 \begin{proof}
 
 By corollary 2.4, $(M,g)$ has uniform finite blocking and is conjugate point free.  By theorem 4.4, $\pi_1(M)$ has polynomial growth.  By \cite{Le}, compact Riemannian manifolds without conjugate points and with polynomial growth fundamental group are flat.
\end{proof}

\section{Finite Blocking Property and Buildings}

In this section, we provide a proof of Theorem 1, which states that
compact quotients of Euclidean buildings have uniform finite blocking.  Let us
start by recalling some elementary facts about Euclidean buildings,
referring the reader to \cite{Br} for more details.

Let $W \subset\mR^n$ be a compact polyhedron, with all faces forming
angles of the form $\pi/m_{ij}$ for some positive integer $m_{ij}$.
Let $\Lm\subset Isom(\mR^n)$ be the Coxeter group generated by
reflections in the faces of the polyhedron, and observe that the
$\Lm$-orbit of $W$ generates a tessellation of $\mR^n$ by isometric
copies of $W$.  We can label the faces of the copies of $W$ in the
tessellation of $\mR^n$ according to the face of $W$ whose orbit
contains them.  A {\it Euclidean building} is a polyhedral complex
$\tilde X$, equipped with a CAT(0)-metric, having the property that
each top dimensional polyhedron is isometric to $W$ (these will be
called {\it chambers}). In addition, a certain number of axioms are
required to be satisfied. We omit a precise definition of Euclidean buildings, 
contenting ourselves with mentioning the properties we will need.  The reader may refer 
to \cite{Br} for a precise definition, and to \cite{Da} for geometric properties of
these buildings.  The polyhedral complex must also satisfy:
\begin{itemize}
\item each face of the complex $\tilde X$ is labelled with one of
the faces of the polyhedron $W$.
\item given any pair of points $x,y\in \tilde X$, there exists an
isometric, polyhedral, label preserving embedding of the tessellated
$\mR^n$ whose image contains $x$ and $y$.  The image of such an
embedding is called an {\it apartment}.
\item the group $Isom(\tilde X)$ is defined to be the group of label
preserving isometries of $\tilde X$.
\item given any two apartments $\mA_1,\mA_2$ whose intersection is
non-empty, there exists an element $\phi\in Isom(\tilde X)$ which
fixes pointwise $\mA_1\cap \mA_2$, and satisfies $\phi(\mA_1)=\mA_2$
\end{itemize}
We will say that $X$ is a {\it compact quotient} of $\tilde X$
provided it is the quotient of $\tilde X$ by a cocompact subgroup of
$Isom(\tilde X)$, acting fixed point freely.

Note that in a Euclidean building, one has uniqueness of geodesics
joining pairs of points (from the CAT(0) hypothesis).  Furthermore,
we can pick an apartment containing both $x$ and $y$, giving a
totally geodesic $\mR^n$ inside $\tilde X$ containing $x,y$.  Then
the geodesic joining $x$ to $y$ coincides with the straight line
segment from $x$ to $y$ within the apartment.  We will call the
point along the geodesic that is equidistant from $x$ and $y$ the
{\it midpoint} of $x$ and $y$, and denote it by $(x+y)/2$.

Another important point is that both the building $\tilde X$, as
well as the compact quotient $X$ come equipped with a canonical {\it
folding map} to the canonical chamber $W$, given by the labeling.
We will use $\rho$ to denote the canonical folding map, and given a
point $x\in X$ (or in $X/\Gamma$), we define the {\it type} of the
point $x$ to be the point $\rho(x)\in W$.  We now make two
observations:
\begin{itemize}
\item in a compact quotient $X/\Gamma$, there are only finitely many
points of any given type, and
\item given any point $p$ in $X/\Gamma$, every pre-image of $p$ in
the universal cover $X$ has exactly the same type as $p$.
\end{itemize}
The proof of the theorem will make use of the following easy:

\begin{Lem}
Let $\mA$ be any apartment, and $x,y\in W$ a pair of points in the
model chamber.  Define $\mS(x)$, $\mS(y)$ to be the set of points in
$\mA$ of type $x$, $y$ respectively.  Then there exist a finite
collection of points $b_1,\ldots ,b_k \in W$ having the property
that:
$$\{(\bar x + \bar y)/2 \hskip 3pt | \hskip 3pt \bar x \in \mS(x), \bar y
\in \mS(y)\} \subset \bigcup_{i=1}^k \mS(b_i)$$
\end{Lem}

\begin{Prf}
We first observe that the Coxeter group $\Lm$ contains an isomorphic
copy of $\mZ^r$ as a finite index subgroup, where $r=\dim (\mA)$. In
particular, if we denote by $\Lp$ this finite index subgroup, we
note that each of the two sets $\mS(x)$, $\mS(y)$ are the union of
$[\Lm : \Lp]=m$ disjoint copies of $\Lp$-orbits in $\mA$.
Now note that given any two $\Lp$-orbits in $\mA$, the
collection of midpoints of points in the first orbit with points in
the second orbit lie in a {\it finite} collection of $\Lp$-orbits
(in fact, at most $2^r$ such orbits). This is immediate from the
proof of the fact that a flat torus has finite blocking.

This in turn implies that the collection of midpoints of points in
the set $\mS(x)$ and points in the set $\mS(y)$ lie in the union of
at most $2^rm^2$ of the $\Lp$-orbits in $\mA$.  Since each
$\Lp$-orbit lies in a corresponding $\Lm$-orbit, we conclude that
the collection of midpoints lie in the union of a {\it finite}
collection of $\Lm$-orbits.  But two points are in the same
$\Lm$-orbit if and only if they have the same type.  Hence choosing
the points $b_1,\ldots, b_k$ to be the finitely many types (at most
$2^rm^2$ of them), we get the desired containment of sets.
\end{Prf}

We now proceed to prove Theorem 1:

\begin{Prf}[Theorem 1]
Let $X=\tilde X/\Gamma$ be a compact quotient of the Euclidean
building $\tilde X$, and let $W$ denote a model chamber.  Given two
points $x,y$ in the space $X$, we want to exhibit a finite set of
blockers.  Consider the sets $\mP(x),\mP(y)\subset \tilde X$
consisting of all pre-images of the points $x,y$, respectively,
under the covering map $\tilde X\rightarrow X=\tilde X/\Gamma$.  As
we previously remarked, we can make sense of the midpoint of a pair
of points in $\tilde X$.  We now claim that the collection of {\it
midpoints} joining points in $\mP(x)$ to points in $\mP(y)$ have
only a finite number of possible types in $W$.

In order to see this, let us apply the previous lemma to the types
$\rho(x),\rho(y)\in W$.  First note that given an arbitrary pair of
points $\bar x\in \mP(x)$, $\bar y\in \mP(x)$, there exists an
apartment $\mA$ containing $\bar x,\bar y$.  Furthermore, we have
the obvious containments $\mP(x)\cap \mA\subset \mS(\rho(x))$,
$\mP(y)\cap A\subset \mS(\rho(y))$, and hence the midpoint joining
$\bar x$ to $\bar y$ has the property that its type is one of the
finitely many points $b_1,\ldots ,b_k$.

Hence to obtain a finite blocking set, let us consider the
collection of {\it all} points in $X$ whose type is one of
$b_1,\ldots,b_k$.  This yields a finite collection $\mathcal B$ of
points in $X$ having the property that if $\gamma$ is an arbitrary
geodesic joining $x$ to $y$, its midpoint must be one of the points
in $\mathcal B$, completing the proof of the theorem.
\end{Prf}


\begin{thebibliography}{99}

\bibitem[Be]{Be} A. Besse.  Manifolds all of whose geodesics are closed. \textit{A Series of Modern Surveys in Mathematics}, \textbf{93}, Springer-Verlag, 1978.



\bibitem[Ber1]{Ber1} M. Berger.  Blaschke's conjecture for spheres.  Appendix D in Manifolds all of whose geodesics are closed.  \textit{A Series of Modern Surveys in Mathematics}, \textbf{93}, Springer-Verlag, 1978.

\bibitem[Ber2]{Ber2} M. Berger.  Lectures on geodesics in Riemannian geometry.  Bombay: Tata Fundamental Institute of Research, 1965.


\bibitem[Br]{Br} K. Brown. Buildings, Springer-Verlag, 1998.


\bibitem[BuGu]{BuGu} K. Burns and E. Gutkin.  Growth of the number of geodesics between points and insecurity for Riemannian manifolds.  Preprint.

\bibitem[Da]{Da} M. Davis.  Buildings are CAT(0).  \textit{Geometry and cohomology in group theory
(Durham 1994)}, Cambridge Univ. Press, 1998, 108-123.

\bibitem[do Ca]{do Ca}  M. do Carmo. Riemannian Geometry, Birkhauser, 1992.

\bibitem[Se]{Se} J.P. Serre. Homologie singuliere des \'espaces fibr\'es. \textit{Ann. Math.}, \textbf{54}, 1951, 425-505.

\bibitem[Fo]{Fo} D. Fomin.  Zadaqi Leningradskih Matematitcheskih Olimpiad, Leningrad, 1990.

\bibitem[Gr]{Gr} L.W. Green. Auf Wiedersehensflachen.  \textit{Ann. Math}, \textbf{78}, 1963, 289-299. 

\bibitem[Gu1]{Gu1}  E. Gutkin.  Blocking of billiard orbits and security for polygons and flat surfaces. \textit{Geom. and Funct. Anal.}, \textbf{15}, 2005, 83-105.

\bibitem[Gu2]{Gu2} E. Gutkin.  Insecurity for lattice translation surfaces of small genus, with applications to polygonal billiards, preprint IMPA D-006, 2005.

\bibitem[Gu3]{Gu3} E. Gutkin.  Blocking of orbits and the phenomenon of (in)security for the billiard in polygons and flat surfaces.  Preprint, IHES/M/03/06.

\bibitem[GuSc]{GuSc} E. Gutkin and V. Schroeder.  Connecting geodesics and security of configurations in compact locally symmetric spaces, \textit{Geometriae dedicata}, \textbf{118} no. 1, 2006, 185-208.

\bibitem[HiSn]{HiSn} P. Heimer and V. Snurnikov. Polygonal billiards with small obstacles. \textit{J. Statist. Phys.}, \textbf{90}, (1998), 453-466.

\bibitem[Le]{Le} N.D. Lebedeva.  On spaces of polynomial growth with no conjugate points. \textit{St. Petersburg Math. J.}, \textbf{16} no. 2, 2005, 341-348.

\bibitem[Ma]{Ma} R. Ma\~n\'e.  On the topological entropy of geodesic flows, \textit{J. Differential Geometry}, \textbf{45}, 1997, 74-93. 

\bibitem[Mo1]{Mo1} T. Monteil. A counter-example to the theorem of Hiemer and Snurnikov., \textit{J. Statist Phys.}, \textbf{114}, 2004, 1619-1623.

\bibitem[Mo2]{Mo2} T. Monteil. On the finite blocking property.  \textit{Annales de l'Institut Fourier},  \textbf{55} no. 4, 2005, 1195-1217.

\bibitem[Mo3]{Mo3} T. Monteil. Finite blocking versus pure periodicity. Preprint, arXiv math.DS/0406506.

\bibitem[Mo4]{Mo4} T. Monteil.  A homological condition for a dynamical and illuminatory classification of torus branched coverings. Preprint, arXiv math.DS/0603352.  

\bibitem[NaRo]{NaRo} A. Nabutovsky and R. Rotman.  The length of the shortest closed geodesic on a 2-dimensional sphere. \textit{Int. Math. Res. Not.}, \textbf{23}, 2002,1211-1222.

\bibitem[Pe]{Pe} J. Pesin.  Formulas for the entropy of the geodesic flow on a compact Riemannian manifold without conjugate points. \textit{Math. Notes}, \textbf{24} no. 4, 1978, 796-805.

\bibitem[Sab]{Sab} S. Sabourau.  Filling radius and short closed geodesics of the 2-sphere. \textit{Bull. Soc. Math. France}, \textbf{132}, no. 1., (2004), 105-136.


\bibitem[Sa]{Sa} L.A. Santalo.  Integral geometry in general spaces. \textit{Proc. Intern. Congress Math. (Cambridge 1950)}, Vol. I., Am. Math. Soc., 1952,483-489. 

\bibitem[Wa1]{Wa1} F. Warner.  The conjugate locus of a Riemannian manifold.  \textit{American Journal of Mathematics}, \textbf{87} no. 3, 1965, 575-604.

\bibitem[Wa2]{Wa2} F. Warner.  Conjugate loci of constant order.  \textit{Annals of Mathematics}, 2nd Ser., \textbf{86} no. 1, 1967, 192-212.

\bibitem[Zo]{Zo} O. Zoll. Uber Flachen mit Scharen geschlossener geodatischer Linien. \textit{Math. Ann.}, \textbf{57}, 1903, 108-133. 

\end{thebibliography}
\end{document}